\documentclass[11pt]{amsart}
\usepackage{amsopn}
\usepackage{amssymb, amscd}

\usepackage{graphics,epsfig}
\newcommand{\nc}{\newcommand}

\nc{\noi}{\noindent}
\nc{\ft}{\footnote}
\nc{\G}{\Gamma}
\nc{\g}{\gamma}
\nc{\Ld}{\Lambda}
\nc{\ld}{\lambda}
\nc{\la}{\langle}
\nc{\ra}{\rangle}
\nc{\ba}{\backslash}
\nc{\ke}{\hspace{-.3cm}}

\renewcommand{\Im}{{\rm Im}}

\nc{\SO}{\mathrm{SO}}
\nc{\Spe}{\mathrm{Sp}}
\nc{\Sl}{\mathrm{SL}}
\nc{\SU}{\mathrm{SU}}
\nc{\Or}{{\mathrm O}}
\nc{\Gl}{\mathrm{GL}}

\nc{\R}{{\mathbb R}}
\nc{\HH}{{\mathbb H}}
\nc{\C}{{\mathbb C}}
\nc{\Z}{{\mathbb Z}}
\nc{\F}{{\mathbb F}}
\nc{\N}{{\mathbb N}}
\nc{\Q}{{\mathbb Q}}
\nc{\PP}{{\mathbb P}}

\nc{\rank}{\operatorname{rank}}

\newtheorem{theorem}{Theorem}[section]
\newtheorem{prop}[theorem]{Proposition}
\newtheorem{cor}[theorem]{Corollary}
\newtheorem{lemma}[theorem]{Lemma}
\newtheorem{definition}[theorem]{Definition}
\newtheorem{remark}[theorem]{Remark}

\title{Generalized Hantzsche-Wendt flat manifolds}

\author[Juan Pablo Rossetti and Andrzej Szczepa\'nski]{J. P. Rossetti$^{(*)}$ and A.
Szczepa\'nski}

\address{FaMAF(Ciem), Universidad Nacional de C\'ordoba, 5000 - C\'ordoba, Argentina.}
\email{rossetti@mate.uncor.edu}
\address{Institute of Mathematics, University of Gda\'nsk,
ul.Wita Stwosza 57, 80 - 952 Gda\'nsk, Poland}
\email{matas@paula.univ.gda.pl}

\thanks{2000 {\it Mathematics Subject Classification.} Primary: 20H15, 57S30;
Secondary: 53C29, 20F34, 57N16, 05C25. \\
{\it Key words and phrases.} Flat manifold, Bieberbach group, holonomy
representation.
\\  $(*)$ Supported by a Guggenheim Fellowship and Conicet}

\begin{document}

\begin{abstract}
We study the family of closed Riemannian $n$-manifolds with holonomy group
isomorphic
to $\Z_2^{\,n-1}$, which we call generalized Hantzsche-Wendt manifolds.
We prove results on their structure, compute some invariants, and find relations
between
them, illustrated in a graph connecting the family.
\end{abstract}

\maketitle

\section{Introduction}

A flat manifold is a closed Riemannian manifold with zero
sectional curvature. From Bieberbach theorems, it is known that in
each dimension there is only a finite number of such manifolds, up
to affine equivalence, and the problem of classifying them is
important. Recently, this has been achieved up to dimension six
(cf.\cite{CS}), but it is not possible to do it in general. Hence
some interesting families of flat manifolds have been considered,
as for instance those with first Betti number zero or those with
holonomy group $\Z_{2}\oplus \Z_{2}$.

The purpose of this article is to study another special class of
flat manifolds, namely those of dimension $n$ with holonomy group
isomorphic to $\Z_{2}^{n-1}$. In dimension two the Klein bottle
belongs to this family and in dimension three there are three such
manifolds: a classical flat manifold first considered by Hantzsche
and Wendt \cite{HW} (now called {\it didicosm}, see last section)
and two non-orientable ones.

Though here they are considered for first time together, these
manifolds have been studied partially in \cite{LS}, \cite{Ma},
\cite{MRc}, \cite{MRb}, \cite{Sz1}, mainly in the oriented case,
where we shall call them {\it Hantzsche-Wendt manifolds} (HW
manifolds for short; respectively HW (Bieberbach) groups, if we
consider their fundamental groups, which are Bieberbach groups).
They have many interesting properties, for instance, there are
pairs of isospectral HW manifolds non homeomorphic to each other
(cf.\ \cite{MRc}). Moreover, HW manifolds have the $\Q$-homology
of spheres (cf.\ \cite{Sz1}) and hence they are geometrically
formal (cf.\ \cite{Ko}).

We shall call the manifolds in our family {\it generalized
Hantzsche-Wendt} manifolds (GHW manifolds for short; respectively
GHW (Bieberbach) groups).

One of the main results in this article is that the integral holonomy
representation of a GHW manifold is of diagonal type (Theorem \ref{diag}).
This simplifies the study of this family.
An immediate consequence is that
the HW manifolds are exactly those considered in \cite{MRc}.

Another interesting property of these manifolds, or their
fundamental groups, is given by the relations between them in
different dimensions. In fact we show that every GHW group has a
subgroup isomorphic to a GHW group of one lower dimension
(Corollary 3.4). Moreover, from Theorem \ref{exist}, any
$n$-dimensional GHW group is a subgroup of some
$(n+1)$-dimensional one. With these observations we define a
graph, whose vertices correspond to GHW groups, connecting the
family, which we draw up to dimension four in the last section.

It is clear that the number of non-oriented GHW flat manifolds is
much greater than the number of oriented ones (see for instance
the table in Section~\ref{no}). In the non-oriented case, there
are two types of GHW groups, depending on whether the first Betti
number equals one or zero. We show that the case with first Betti
number one (as well as the oriented case) has only one holonomy
representation, which is completely determined. Hence it makes
sense to carry out the calculation of the $\Q$-homology (see
Theorem~\ref{betti}).

The situation in the non-oriented case with first Betti number
zero is much more complicated since the holonomy representation is
not unique. We find the exact number of integral representations
occurring in this subclass and show that it grows linearly with
the dimension. This follows as a consequence of Theorem
\ref{exist}, which asserts that every diagonal integral
representation can be realized as the integral holonomy
representations of a flat manifold. In Section \ref{examples} we
consider some examples in this subclass which are similar to the
HW manifolds from \cite{Sz1}. In Section \ref{outaut} we study
some symmetries of GHW manifolds. In particular, we determine the
outer automorphism groups of certain subfamily of GHW groups with
non-trivial center. Moreover a general upper bound of the order of
the outer automorphism group of GHW groups is given.

There are some problems and questions which we formulate at the end.

We would like to thank Wim Malfait for discussions regarding
Section~\ref{outaut}.

\section{Preliminaries and Notations}\label{pn}

Let $\Gamma$ be a Bieberbach group of dimension or rank $n$. By definition it is a
torsion-free, discrete, cocompact subgroup of the group $\mathrm{E}(n)=
\Or(n)\ltimes\R^{n}$ of rigid motions
of the $n$-dimensional Euclidean space.
The multiplication in $\mathrm{E}(n)$ is given by
$$(B,b)(C,c) = (BC,Bc+b),$$
where $B,C\in \Or(n)$ and $b,c\in\R^n$. The translation part $c$
of an isometry $(C,c)\in \mathrm{E}(n)$ shall sometimes be denoted
by $s(C)$. If a flat manifold $\R^{n}/\Gamma$ is oriented then by
definition $\Gamma$ is oriented. From Bieberbach theorems (cf.\
\cite{Ch}), the translation subgroup $\Ld\simeq\Z^{n}$ of $\Gamma$
is normal, maximal abelian, of finite index. Then there is a short
exact sequence of groups
\begin{equation}\label{short}
0 \longrightarrow \Ld \longrightarrow \Gamma \longrightarrow G \longrightarrow
0,
\end{equation}
where $G$ is finite and isomorphic to the holonomy group of the
manifold $\R^{n}/\Gamma$, thus we call it the holonomy group of
$\Gamma.$ It is well known that any extension as (\ref{short})
defines an element $\alpha \in H^{2}(G,\Ld)$, where the $G$-module
structure of $\Ld$ is defined by conjugation in (\ref{short}). We
denote by $\rho$ this integral representation and call it the {\it
integral holonomy representation} of $\G$ (or of $\R^n/\G$).

\begin{definition}{\rm
Let $\Gamma$ be a Bieberbach group of dimension $n$.
A flat manifold $\R^{n}/\Gamma$ is called {\it generalized Hantzsche-Wendt}
manifold
if its holonomy group is $\Z_2^{\,n-1}$.
In this case $\G$ is called a {\it generalized Hantzsche-Wendt} group.
}\end{definition}

\begin{definition}{\rm
We say that a Bieberbach group $\Gamma$ is {\it diagonal} or of {\it diagonal
type} if its
corresponding lattice of translations $\Ld$ has an orthogonal basis
in which the integral holonomy representation $\rho$ diagonalizes.
We will also call {\it diagonal} the corresponding flat manifold
$\R^{n}/\Gamma$.}
\end{definition}

\begin{remark}{\rm
In \cite{MRd}, the previous definition was done with the (slightly
more restrictive) assumption that the basis of $\Ld$ is
orthonormal. We recall that diagonal manifolds have necessarily
the holonomy group isomorphic to $\Z_2^{\,k}$ for some $k\in\N$.
}\end{remark} Consider
\begin{equation}\label{Ci}
C_i := \left[ \begin{array}{cccccccc}
1 & 0 & \cdots & 0 & 0 & 0 & \cdots & 0 \\
0  & 1 &  \cdots & 0 & 0 & 0 & \cdots & 0  \\
\vdots & \vdots & \ddots &  &  \vdots & \vdots &  &
\vdots \\
0 & 0 &  & 1 & 0 & 0 & \cdots & 0  \\
0 & 0 & \cdots & 0 & \!\!-1 & 0 & \cdots & 0  \\
0 & 0 & \cdots & 0 & 0 & 1 &  & 0  \\
\vdots &  \vdots & &  \vdots  & \vdots & & \ddots &
\vdots \\
0 & 0 & \cdots & 0 & 0 & 0 & \cdots & 1
\end{array}\right] , \,\,\, 1\leq i \leq n,
\end{equation}
where the $-1$ is placed in the $(i,i)$ entry.
For brevity we will denote
$C_i=\hbox{diag}(1,\dots,1,\underbrace{-1}_i,1,\dots,1)\,$, for $\,1\le i\le
n$.
The product of two of these is
\begin{equation}\label{2prod}
C_iC_j= {\rm
diag}(1,\dots,1,\underbrace{-1}_i,1,\dots,1,\underbrace{-1}_j,1,\dots,1),\,
\hbox{
for every } i<j.
\end{equation}
The elements in (\ref{2prod}) generate the group $\mathrm{D}^+(n)
:= \mathrm{D}(n) \cap {\Sl}(n,\Z)$, of cardinality $2^{n-1}$, and
of index two in the group $\mathrm{D}(n)$ of diagonal $n\times n$
matrices with entries $\pm 1$ in the diagonal.

In the oriented case with holonomy group $\Z_2^{n-1}$, the
rationalization of the integral holonomy representation $\rho\,$
(i.e.\ to consider $\Ld$ as a $\Q$-module instead of a
$\Z$-module), gives a unique representation $\rho_\Q$ whose image,
in an appropriate basis, is exactly $\mathrm{D}^+(n)$. Hence, when
$n$ is even there are no HW manifolds, since $\mathrm{D}^+(n)$
contains $-Id$, which always corresponds to an element of torsion.

\section{Diagonal Structure of GHW Manifolds}

To start with, we present an important structural property of the class of
GHW
manifolds.

\begin{theorem}\label{diag}
The fundamental group of a Generalized Hantzsche-Wendt manifold is diagonal.
\end{theorem}

For the proof we start with two lemmas.

\begin{lemma}\label{contains}
Let $\rho:\Z_2^{\,n-1}\rightarrow \Gl(n,\Z)$ be a diagonal faithful integral
representation with $-Id\not\in\Im(\rho)$.
Then there is $g\in\Z_2^{\,n-1}$ such that
$\rho(g)=\mathrm{diag}(-1,\dots,-1,1,-1,\dots,-1)$.
Moreover, if furthermore $\,\Im(\rho)\not\subseteq \Sl(n,\Z)$ then there is
$g\in\Z_2^{\,n-1}$ such that  $\rho(g)=\mathrm{diag}(1,\dots,1,-1,1,\dots,1)$.
\end{lemma}
\begin{proof}
First, we regard $\mathrm{D}(n)$, $\mathrm{D}^+(n)$ and
$\Im(\rho)$ as $\Z_2$-vector spaces, of dimensions $n$, $n-1$ and
$n-1$ respectively. Also, $\mathrm{D}(n) = \Im(\rho) \cup
(-Id)\Im(\rho)$. Since the $C_i$'s, for $1\le i\le n$, are
linearly independent, they cannot lie all in $\Im(\rho)$
simultaneously, thus there is at least one of them in
$(-Id)\Im(\rho)$, or equivalently, there is a $-C_i$ in
$\Im(\rho)$, as claimed.

The last assertion in the statement follows from a similar
argument.
\end{proof}

\begin{lemma}\label{split}
Let $\Gamma$ be a $n$-rank Bieberbach group with translation lattice $\Ld$.
Suppose that $(B,b)\in\Gamma$ and $B$ has eigenvalues $1$, $-1$, with
corresponding eigenspaces $V^+$ and $V^-$ of dimensions $1$ and $n-1$
respectively.
Then $\Ld=(\Ld\cap V^+)\oplus (\Ld\cap V^-)$, and the orthogonal projection of
$b$
onto $V^+$ lies in $\frac 12 (\Ld\cap V^+)\smallsetminus (\Ld\cap V^+)$
\end{lemma}
\begin{proof}
We have
\begin{equation}\label{square}
\G\ni {(B,b+\ld)}^2=\left(Id,(B+Id)(b+\ld)\right), \textrm{ for every }
\ld\in\Ld.
\end{equation}
The torsion-free condition implies that $0\neq
(B+Id)(b+\ld)\in\Ld$, for every $\ld\in\Ld$. Now, $\frac 12(B\pm
Id)$ is the orthogonal projection onto $V^{\pm}$. If we suppose
that for some $\ld\in\Ld$, $\frac12 (B+Id)(b+\ld)=\mu\in\Ld$, then
$(B+Id)(b+\ld-\mu)= (B+Id)(b+\ld)- (B+Id)(\mu)=2\mu-2\mu=0$, a
contradiction. Hence, the orthogonal projection of $b$ (and also
of $b+\ld$ for every $\ld\in\Ld$) onto $V^+$ lies in $\frac 12
\Ld\smallsetminus\Ld$. The fact that $V^+$ is one-dimensional now
implies that the orthogonal projection of $\Ld$ onto $V^+$ is
exactly $\Ld\cap V^+$. Hence $\Ld$ splits into the direct sum of
two lattices, and the lemma follows.
\end{proof}

\smallskip
\noindent {\it Proof of Theorem \ref{diag}.} We proceed by
induction on the dimension $n$. If $n$ is small, it is known that
the theorem holds. Let $\rho$ be the integral representation of
$\G$ defined on $G=\Z_2^{\,n-1}$ and let $\Ld \simeq \Z^{n}$ be
the lattice of translations of $\Gamma$. The nature of the
holonomy group $G$ implies that each $\rho_{\Q}(g)$ diagonalizes,
and since they commute, they diagonalize simultaneously, i.e.,
there is a $\Q$-basis $\mathcal B$ of $\Ld$ such that
$\left[\rho_{\Q}(g)\right]_{\mathcal B}$ is diagonal for every
$g\in G$ (and the elements in the diagonal are $\pm1$'s). Then,
Lemma \ref{contains} applied to these integral matrices implies
that $\exists \, g_0\in\Z_2^{\,n-1}$ such that
$\left[\rho_{\Q}(g_0)\right]_{\mathcal
B}=\textrm{diag}(1,-1,\dots,-1)$ (after reordering $\mathcal B$).
Hence, $\rho(g_0)$ has eigenvalues $+1$ and $-1$ with
multiplicities $1$ and $n-1$ respectively. Now we can a apply Lemma
\ref{split}, and thus $\Ld$ splits into a direct sum
$\Ld=\Ld_1\oplus\Ld_2$ where dim$\Ld_1=1$ and dim$\Ld_2=n-1$. Also
$\textrm{span}_\R \Ld_1$ (resp.\ $\textrm{span}_\R \Ld_2$) is the
eigenspace of $\rho(g_0)$ with eigenvalue $1$ (resp.\ $-1$) and
the orthogonal projection of the translational part corresponding
to $\rho(g_0)$ lies in $\frac 12 \Ld_1\smallsetminus\Ld_1$. Since
$\rho(g_0)\rho(g) = \rho(g)\rho(g_0)$ for every $g\in
\Z_2^{\,n-1}$, we have that $\Ld_1$ and $\Ld_2$ are stable under
the action of $\rho(g)$, thus $\Ld=\Ld_1\oplus\Ld_2$ is a
$\Z_2^{\,n-1}$-module decomposition. Hence we have that
$\rho=\chi\oplus\rho'$, where $\chi$ is a character and $\rho'$
has rank $n-1$, both defined on the same domain as $\rho$. Thus,
$\forall g\in\Z_2^{\,n-1}$

\

\begin{tabular}{r|c|c|  c c c  r|c|c| }
\cline{2-3}\cline{8-9}
 & $\chi(g)\!$ &    & & & &  & \,\,1\,\, & \\
\cline{2-3} \cline{8-9}
$\rho(g)=$ & &    & & \quad and \quad &  & $\rho(g_0)=$
& & \\
 & & $\rho'(g)$ & & & & & & $-Id$  \,\\
 & &    & & &   & & & \\
\cline{2-3}\cline{8-9}
\end{tabular}

\

We may conjugate $\Gamma$ by $(Id,\mu)$ with $\mu\in \frac 14\Ld_1$ and obtain
that
$p_1(b)\in \frac 12\Ld_1$ for every $(B,b)\in\Gamma$,
where $p_1$ is the orthogonal projection onto span$_\R\Ld_1$ (cf.\ \cite{MRd},
Lemma 1.4).

We define a character $\psi:{\Gamma/\Ld} \simeq\Z_2^{\,n-1}
\longrightarrow\Z_2=\{0,1\}\,$ by

\begin{displaymath}
\psi((B,b))= \left\{
\begin{array}{ll}
0 & \textrm{ if } \,\,p_1(b)\in\Ld_1; \\
1 & \textrm{ if } \,\,p_1(b)\in\frac 12\Ld_1\smallsetminus\Ld_1.
\end{array} \right.
\end{displaymath}
\noindent It is well defined since all the elements $(B,b)$ in the same class
in $\Gamma/\Ld$ have the same $p_1(b)$ mod $\Ld_1$.
Now we consider
$$\Gamma \stackrel{r}{\longrightarrow} \Gamma/\Ld
\stackrel{\psi}{\longrightarrow} \Z_{2},$$ where $r$ is the
canonical projection. Set $\Gamma':=\ker(\psi\circ r)$. It is a
subgroup of $\Gamma$ of index~2 and its elements $(B,b)$ have the
property that $p_1(b)\equiv 0 \,\, (\textrm{mod }\Ld_1)$.  The
last fact shows that if we take $\Gamma'':=\Gamma'\mid_{\Ld_2}$
(i.e. the restriction to the last $(n-1)$ coordinates of $B$'s and
$b$'s), it is still a Bieberbach group, it has rank $n-1$ and
holonomy group isomorphic to $\Z_{2}^{\,n-2}$. (Alternatively, we
could have considered the subgroup $\Gamma_{0}:=\{(B,b)\in\Gamma :
b_{1}=0\}$, where $b=(b_{1},b_{2})\in\Lambda_1\oplus\Lambda_2$,
and observed that the restriction  of $\Gamma_0$ to
$\textrm{span}_\R \Ld_2$ gives an isomorphism between $\Gamma_{0}$
and $\Gamma''$.) Then we can apply induction to $\Gamma''$, hence
it is diagonal. By the relation between $\Gamma$ and $\Gamma''$ it
is clear that therefore $\Gamma$ is diagonal. \qed
\smallskip

In the proof, the group denoted by $\G''$ is a subgroup of $\G$
(in the abstract sense, i.e., there is a monomorphism from $\G''$
into $\G$). This gives a way of going from one dimension to the
previous one within the class of GHW manifolds, which will be used
in Section~\ref{graph}. In other words we have
\begin{cor}\label{subgroup}
Every GHW group has a subgroup isomorphic to a GHW group of one lower rank.
\end{cor}

In the oriented case, the existence of these manifolds in each odd
dimension was proved in \cite{Ma} and one realization in each odd
dimension was shown in \cite{Sz1}. In \cite{MRc} the oriented flat
manifolds with (orthonormal) diagonal holonomy and holonomy group
$\Z_2^{n-1}$ were studied and they were parametrized by certain
directed graphs; also lower bounds for their cardinality were
shown, growing rapidly with $n$.

\begin{remark}\label{iii}{\rm
(i) A direct consequence of Theorem~\ref{diag} is that in the
oriented case, the above definition of generalized Hantzsche-Wendt
manifolds coincides with the Hantzsche-Wendt manifolds considered
in \cite{MRc}.

(ii) Generalized Hantzsche-Wendt manifolds are of minimal
dimension among those with holonomy group $\Z_2^{\,k}, k$ fixed
and $\geq 1$ (cf.\ \cite{HMSS}). Furthermore they are the only
existing closed Riemannian $n$-manifolds with holonomy group
$\Z_2^{\,n-1}$, since manifolds with finite holonomy group must be
flat, by a theorem of Cartan-Ambrose-Singer (see \cite{Wo}, Cor.\
3.4.7).

(iii) For every $k$, $1<k<n-1$, there are flat manifolds
with holonomy group $\Z_2^{\,k}$ not of diagonal type. Examples of this
situation are
easy to construct (see end of Section \ref{examples}).

(iv) Somehow opposite to diagonal manifolds are flat manifolds
having indecomposable holonomy representation. We recall that
there are many of them, and for example, for holonomy group
$\Z_2^{\,2}$, they have been classified and the set of their
dimensions is unbounded (cf.\ \cite{BGR}). }\end{remark}

\section{Non-Oriented GHW Flat Manifolds}\label{no}

The class of non-oriented GHW manifolds
includes two subclasses, depending on whether the first Betti number $\beta_1=0$
or~1,
or equivalently (cf.\ \cite{HaS}), the classes of Bieberbach groups with
trivial or non-trivial center.

Let us consider first the case with non-trivial center. From
Theorem \ref{diag} and a similar observation as at the end of
Section~\ref{pn} for HW manifolds, this subclass have the integral
holonomy representation uniquely determined. A simple example is
the following (cf.\ \cite{LS}). Let $K_n$ be the subgroup of
$\mathrm{E}(n)$ generated by the set \,\{$(C_{i}, s(C_{i})) :
0\leq i \leq n-1\}$, where $C_{0} = I$ and
\,$C_i=\text{diag}(1,\dots,1,\underbrace{-1}_i,1,\dots,1)\,$ are
the $n\times n$ orthogonal matrices defined in (\ref{Ci}).
Moreover \,$s(C_{0}) = e_{1}$ and \,$s(C_{i}) = e_{i+1}/2$ for
other $i$. Here $e_{1},e_{2}, ... ,e_{n}$ denote the canonical
basis of \,$\Z^{n}$. In \cite{LS} it was proved that $K_n$ is a
Bieberbach group and that there exists the following short exact
sequence of groups
\begin{equation}\label{lee}
0\longrightarrow K_{n-1}\longrightarrow K_n\longrightarrow \Z \longrightarrow
0,
\end{equation}
where the epimorphism $\,K_n\rightarrow \Z\,$ is given by the composition
$$K_n\rightarrow K_n/[K_n,K_n]\rightarrow\Z.$$
However in general we have the following proposition.

\begin{prop}\label{chain}
Let $\Gamma$ be a Generalized Hantzsche-Wendt group of dimension $n$
with non-trivial center. Then the kernel of the composition
$f:\Gamma \rightarrow \Gamma/[\Gamma,\Gamma]\rightarrow \Z$ can be an
oriented or
a non-oriented Generalized Hantzsche-Wendt group.
\end{prop}
\begin{proof} First, it is not difficult to observe that $\ker f$ is a
Bieberbach group (cf.\ \cite{Ch}). From the assumptions, the
abelian group $\Gamma/[\Gamma,\Gamma]\simeq \Z \oplus F$ where $F$
is a finite group (cf.\ \cite{HaS}). Let $f((A,a))\in \Z$ be a
generator for some $(A,a)\in \Gamma$. We can assume, by using
Lemma~\ref{contains}, that $A$ has eigenvalues $1, -1$, with
corresponding eigenspaces $V^{+}$ and $V^{-}$ of dimensions $1$
and $n-1$ respectively. Moreover we have an isomorphism $\ker
f/(\Lambda \cap \ker f) \simeq \ker f/(\Lambda \cap V^{-}) \simeq
\Z_2^{\,n-2}$, where $\Lambda$ is the translation lattice of
$\Gamma$. We have to show that $\Z_2^{\,n-2}$ is the holonomy
group of the Bieberbach group $\ker f$. But this follows from the
assumptions and Lemma \ref{split}. To finish the proof we observe
that for odd $n$, every $\ker f$ is non-oriented. On the other
hand, using semidirect product $\G'\rtimes_{\alpha}\Z$ with the
integers, for any HW group $\G'$ of dimension $n$, we can define a
$(n+1)$-rank GHW group with non-trivial center such that the
kernel of the map onto $\Z$ is isomorphic to~$\G'$. Here $\alpha :
\Z \rightarrow \mathrm{Out}(\G')$ is a homomorphism which send $1$
to $-Id \in \mathrm{Out}(\G').$
\end{proof}

As mentioned above, the holonomy representation of GHW manifolds with
$\beta_1=1$ is always
equivalent to that of the above example, hence all these manifolds have the
same rational homology,
namely

\begin{theorem}\label{betti}
Let $M$ be a Generalized Hantzsche-Wendt flat manifold of dimension $n$
with first Betti number one. Then
$$H_{j}(M,\Q) \simeq \left\{ \begin{array}{ll}
\Q \,\,\,& \text{if } \,j=0,\,1; \\
0 & \text{otherwise.}\\
\end{array} \right.$$
\end{theorem}

\begin{proof} From the definition and assumptions we have a short exact
sequence
$$
0\longrightarrow \Z^{n}\longrightarrow \pi_{1}(M)\longrightarrow
\Z_2^{\,n-1}\longrightarrow 0,
$$
where the action of \,$\Z_2^{\,n-1}$ on $\Z^{n}$ is given by the matrices
$C_{i}$,
$1\leq i \leq n-1.$ This defines a representation
$\rho:\Z_2^{\,n-1}\rightarrow \Gl(n,\Z).$ Using a similar method as in
\cite{Sz1}
we
have a character $\Phi$ of $\rho$ given by the formula
$$\Phi = \sum_{i=1}^{n} \chi_{i}\,,$$
where $\chi_{i} = \chi'_{0}\dots \chi'_{0}\chi'_{1}\chi'_{0}\dots
\chi'_{0}, 1\leq i \leq n-1$ with $\chi'_{1}$ the $i$-th of the
$(n-1)$ factors in the product, and $\chi_{n} =
(\chi'_{0})^{n-1}.$ Here $\chi'_{0}, \chi'_{1}$ denote
respectively a trivial and a non-trivial character of the
involution group $\Z_{2}$. From \cite{Sz1}
$$H_{j}(M,\Q) \simeq [\Lambda^{j}(\Q^{n})]^{\Z_2^{\,n-1}}.$$
From its definition, the representation of $\Z_2^{\,n-1}$ on
$\Lambda^{j}(\Q^{n})$ has the following character
\begin{equation}\label{lambda}
\Lambda^{j}(\rho) = \sum_{1\leq i_{1}< i_{2} < \dots < i_{j} \leq n}
\chi_{i_{1}}\chi_{i_{2}}\dots \chi_{i_{j}}.
\end{equation}
We have
$$\Lambda^{j}(\rho) = \chi_{1}\sum_{2\leq i_{2} < \dots < i_{j} \leq n}
\chi_{i_{2}}\dots \chi_{i_{j}}
+ \sum_{2\leq i_{1}< i_{2} < \dots < i_{j} \leq n}
\chi_{i_{1}}\chi_{i_{2}}\dots \chi_{i_{j}}.
$$
Hence by induction and an immediate calculation the result follows.
\end{proof}

According to \cite{CS} the number of isomorphism classes of GHW manifolds
and their holonomy representations in dimensions up to 6 are as in the
following
table.

\

\begin{center}
\begin{tabular}{|c|r|r||r||r|r|c|}
\hline
 & \multicolumn{5}{|c|} {number of GHW manifolds.} & {number of}     \\
\cline{2-6}
dim & $\beta_1=0$ & $\beta_1=1$ &  total & orient. &
non-orient. &  holon.\ repr. \\
\hline
2 & 0 & 1 & 1 & 0 & 1 & 1 \\
3 & 1 & 2 &  3 & 1 & 2 & 2\\
4 & 2 & 10 & 12 & 0 & 12 & 2 \\
5 & 23 & 100 & 123 & 2 & 121 & 3 \\
6 & 352 & 2184 & 2536 & 0 & 2536 & 3 \\
\hline
\end{tabular}
\end{center}

\

There are 62 HW manifolds in dimension 7 and more that 1000 in
dimension 9 (cf.\ \cite{MRb}).

\section{Holonomy Representations}

If we consider the case of non-oriented generalized Hantzsche-Wendt
Bieberbach groups with trivial center then we have many more possibilities
for the holonomy representation.
Let us start establishing the following general result.

\begin{theorem}\label{exist}
Let $\rho:\Z_2^{\,k} \rightarrow \Gl(n,\Z)$ be a diagonal faithful integral
representation with $-Id\not\in\Im(\rho)$ and $k<n$.
Then there is a Bieberbach group whose integral holonomy representation is
$\rho$.
\end{theorem}

\begin{proof}
It suffices to prove the case $k=n-1$, since in the cases $k<n-1$ we can
use
the lemma below to induce from $\rho$ a new
representation defined on $\Z_2^{\,n-1}$; then apply the
result for $n-1$ to obtain a diagonal Bieberbach group $\Gamma$ with the new
integral holonomy representation;
and finally, choose the appropriate proper subgroup of $\Gamma$ of index
$2^{n-1-k}$ corresponding to
the original integral representation $\rho$.
Note that this chosen subgroup has the same lattice of translations as
$\Gamma$,
and by construction, it is always a diagonal Bieberbach group.

\begin{lemma}
If $\sigma:\Z_2^{\,k} \rightarrow \Gl(n,\Z)$ is a diagonal faithful integral
representation with $-Id\not\in\Im(\sigma)$
and $k<n-1$, then there exists an integral representation
$\tau:\Z_2^{\,k+1} \rightarrow \Gl(n,\Z)$ with the same
properties as $\sigma$ and containing properly $\sigma$, i.e.,
$\tau\restriction_{\Z_2^{\,k}\times\{0\}}=\sigma$.
\end{lemma}
\begin{proof}
Denote $\sigma^-:=-Id\circ\sigma$, i.e.
$\sigma^-(g):=-(\sigma(g))$. Since $-Id\not\in\Im(\sigma)$, we
have $\Im(\sigma)\cap\Im(\sigma^-)=\emptyset$. Since
$2\times2^k<2^n$, there is $B\in\mathrm{D}(n)\smallsetminus
(\Im(\sigma)\cap\Im(\sigma^-))$. Now define $\tau:\Z_2^{k+1}\to
\mathrm{D}(n)$ as the unique representation satisfying
$\tau\restriction_{\Z_2^k\times\{0\}}=\sigma$ and
$\tau\restriction_{\{0\}^k\times\Z_2}=\langle B\rangle$ (here we
are identifying $\langle B\rangle$ with the representation of
$\Z_2$ in $\Z^n$ whose nontrivial element acts by $B$). Since
$B\not\in\Im(\sigma)$, $\tau$ is a faithful representation; since
$B\not\in\Im(\sigma^-)$, $-Id\not\in\Im(\tau)$, and therefore
$\tau$ has the desired properties.
\end{proof}

\vskip.02cm
\noindent{(\it continuation of the proof of Theorem \ref{exist})}
 We continue the proof by induction on $n$.
If $n=2,3$, it is easy to prove it. If $\Im(\rho)\subseteq
\Sl(n,\Z)$, then $n$ must be odd and we know that there are
manifolds with $\rho$ as its integral representation, namely the
Hantzsche-Went manifolds as in \cite{Ma}, \cite{MRc}, \cite{Sz1}.
Otherwise $\Im(\rho)$ contains elements of determinant $-1$, then
by Lemma \ref{contains}, there is $C_i\in \Im(\rho)$, that after
reordering the coordinates we may assume to be
$C_1=\textrm{diag}(-1,1,\dots,1)$. Denote $p_1$ the projection
onto the first coordinate and $p_r$ the projection onto the last
$n-1$ coordinates. Next, consider $\ker(p_1\circ\rho)$, which is a
subgroup of $\Z_2^{\,n-1}$ of index~2. Denote the inclusion
$\ker(p_1\circ\rho)\buildrel{\iota}\over\hookrightarrow
\Z_2^{\,n-1}$, and define $\rho':=p_r\circ\rho\circ\iota$. Then
$\rho'$ is faithful, diagonal and $-Id\not\in\Im(\rho')$, however
it is of rank $n-1$ and defined from a group isomorphic to
$\Z_2^{\,n-2}$. By induction, there is a Bieberbach group
$\Gamma'$ whose integral holonomy representation is $\rho'$.
Moreover, we may assume that all the elements in $\Gamma'$ are of
the form $(B',b')$ with $B'\in \Im(\rho')$ and $b'\in\frac
12\Z^{n-1}$ (see Lemma 1.4 in \cite{MRd}).

To complete the proof one might use previous known results, as in
\cite{Sz2}, but we prefer the following explicit method, since it
is elementary and provides a method which will be useful in
Section~\ref{graph}.

\smallskip

\begin{center}
\begin{tabular}{|c c c |c| c c c|}
\hline
1 & $\cdots$ & 1 & \!$-1$\, & \!$-1$\, & $\cdots$ &
$-1$
\\
\hline
 & & & 1 & & & \\
& $\mathbf{\G'}$ & & \vdots & & $\mathbf{\G'}$ & \\
& & & 1 & & & \\
\hline
\end{tabular}
\end{center}

\smallskip

The previous picture gives an idea of the Bieberbach group we are constructing.
In the middle column is
placed $C_1$.
We need a vector $c\in\R^n$ to form a torsion-free element $(C_1,c)$.
We choose $c={e_2}/2$ (any other choice $c={e_i}/2$ for $i>2$ would work
well
too).
Now we define new elements by using the elements in $\G'$ in the following way:
if $(B',b')\in\G'$, then we define $(B,b)$ by $B:=\left[\begin{array}{cc} 1
& \\  & B'\end{array}\right]$
and $b:=(b_1',0,b_2',\dots,b_{n-1}')$,
where $b'_i$, $1\le i\le n$, are the coordinates of $b'$.
Then $\widetilde{\G'}:=\langle (B,b) :  (B',b')\in\G'; \, L_{e_i}:1\le i\le n
\rangle$
is an $n$-rank Bieberbach group with holonomy group $\ker(p_1\circ\rho)$ (it is
torsion-free since $\G'$ is so).
We claim that $\langle (B,b)\in\widetilde{\G'}, (C_1,c)\rangle$ is a Bieberbach
group with
integral holonomy representation $\rho$.
It is clear that its integral representation is $\rho$, so the only point to be
explained is that it
is torsion-free. This condition is clear for elements in $\widetilde{\G'}$ and
for $(C_1,c)$.
Now the products $(B,b)(C_1,c)$ are of the forms $(D,d)$ with
$D=\left[\begin{array}{cc} -1 & 0 \\ 0 & B' \end{array}\right]$
and $d=(\ast,\frac 12,b'_2,\dots,b'_{n-1})$.
Since $\G'$ is torsion-free, then the elements $(D,d)$ are torsion-free too.
\end{proof}

Now we are interested in the holonomy representations of GHW groups.

\begin{theorem}\label{number}
The total number of inequivalent integral holonomy
representations occurring in Generalized Hantzsche-Wendt Bieberbach groups of
dimension $n$ equals $n/2$ for $n$ even and $(n+1)/2$ for $n$ odd.
\end{theorem}

\begin{proof}
By definition, each holonomy representation is defined by a
hyperplane in $\Z_2^{\,n}$. From linear algebra, it is the kernel
of a linear map with matrix $[x_{1}, x_{2},\dots ,x_{n}]$, where
$x_{i}\in \Z_{2}=\{0,1\}$. Hence, we have the following $n$
classes of diagonal representations:
$$[1,0,\dots,0], \,[1,1,0,\dots,0],\dots,[1,1,\dots,1,0],
\,[1,1,\dots,1].$$ Since we are considering torsion-free groups,
all the hyperplanes corresponding to holonomy representations do
not have the vector $[1,1,\dots,1]$, since $-Id$ corresponds to an
element of torsion. Hence we have to consider only vectors with an
odd number of nonzero elements. We claim that diagonal holonomy
representations defined by such vectors are all inequivalent. In
fact, let us consider a hyperplane (representation) as the set of
$n \times n$ diagonal matrices. It suffices to count the number of
matrices $A$ such that the group of fix points $(\Z^{n})^{A}$ of
the action of $A$ is isomorphic to the integers~$\Z$, i.e.
matrices of the form $A=\mathrm{diag}(-1,\dots,-1,1,-1,\dots,-1)$.
It is not difficult to see that for different hyperplanes, the
number of such $A$'s is different. Summing up we have for $n$
even, $n/2$ and for $n$ odd, $(n+1)/2$ possibilities, which indeed
occur by Theorem \ref{exist}.
\end{proof}

As a consequence, the number of integral representations occurring
in non-oriented GHW groups with trivial center is $(n/2)-1$ for
$n$ even and $((n+1)/2)-2$ for $n$ odd.

\section{Examples}\label{examples}

In this section we present an explicit example of GHW groups with
trivial center in each dimension $n$. In the notation of the proof
of Theorem~\ref{number}, its holonomy representation corresponds
to hyperplanes $[1,1,\dots,1,1]$ in the odd case and to
$[1,1,\dots,1,1,0]$ in the even case. For $n\ge 2$, let
$\Gamma_{n}$ be the subgroup of $\mathrm{E}(n)$ generated by the
set $\{(B_{i},s(B_{i}) = x_i): 1\leq i \leq n-1\}.$ Here $B_{i}$'s
are the $n\times n$ orthogonal matrices:
\begin{equation}\label{Bi}
B_i:=-C_i=\text{diag}(-1,\dots,-1,\underbrace{1}_i,-1,\dots,-1)
\end{equation}
and
\begin{center}
$x_i=s(B_{i}) = e_{i}/2 + e_{i+1}/2$\, for \,$1 \leq i \leq n-1$.
\end{center}
These groups have trivial center for $n\ge 3$.
In \cite{Sz1} (see also \cite{MRc}) it was proved that $\Gamma_{n}$ are
Bieberbach groups for odd $n$.
For even $n$ the groups $\Gamma_{n}$ are also Bieberbach groups
(one might also use Proposition 2.1 in \cite{DM} to prove it)
in particular because we have

\begin{prop}\label{mono}
For $n \geq 2,$ there is a monomorphism $\phi_{n}:\Gamma_{n}\rightarrow
\Gamma_{n+1}$ of groups. Moreover $\phi_{n}(\Gamma_{n})$ is not a normal
subgroup of $\Gamma_{n+1}$.
\end{prop}

\begin{proof}
Let \,\{$(B_{i}, x_{i}) : 1 \leq i \leq n-1\}$\, be a set
of generators of $\Gamma_{n}.$ We have an inclusion $k:\Z^{n}\rightarrow
\Z^{n+1}$, where $b\in \Z^{n}$ is mapped to $k(b) = (b,0)\in \Z^{n+1}.$
Next, we consider the $(n+1)\times(n+1)$ matrices
$\bar{B_{i}}:= \left[ \begin{array}{c|c}
B_i & 0 \\
\hline
0 & -1
\end{array}\right], \,1\leq i \leq n.$
We define
$$\phi_{n}((B_{i}, x_{i})) = (\bar{B_{i}}, k(x_{i})).$$
Finally, it suffices to conjugate by the translation
$(Id,(0,\dots,0,1))$ an element
$\phi_{n}((B_{1},(1/2,1/2,0,\dots,0))$ to show that the image of
$\Gamma_{n}$ is not a normal subgroup.
\end{proof}

By analogy with \cite{Sz1} we have.
\begin{theorem}\label{betti'}
For $n \geq 2$,\,
$$H_{j}(\R^{n}/\Gamma_{n},\Q)\simeq
\left\{ \begin{array}{cll}
\Q &  \text{ if }\,\, j=0 \\
\Q &  \text{ if }\,\, n \text{ is even and } j=n-1\\
\Q &  \text{ if }\,\, n \text{ is odd and } j=n\\
0 &   \text{otherwise}
\end{array} \right.$$
\end{theorem}

\begin{proof}
For odd $n$ it was proved in \cite{Sz1}.
The proof for $n=2k\,$ follows the same structure as the proof of Theorem
\ref{betti}.
We have just to replace $n$ by $2k$, $M$ by $\G_{2k-1}$, and the action of
matrices $C_i$ by the action of $B_i$
as in (\ref{Bi}), thus the characters change accordingly.
Now from the equivalent to formula (\ref{lambda}), an immediate calculation
shows that
no summand is trivial there, then
dim$_{\Q}[\Lambda^{j}(\Q^{2k})]^{\Z_2^{\,2k-1}} = 0$ for each $j$, $1\leq j\leq
2k-2$.
Moreover
$$
\sum_{g\in \Z_2^{\,2k-1}} \frac 1{2^{2k-1}}  \sum_{1\leq
i_{1}<i_{2}<\dots<i_{2k-1}\leq 2k}
\chi_{i_{1}}\chi_{i_{2}}\dots\chi_{i_{j}}(g) =$$
$$=\!\!\sum_{g\in\Z_2^{\,2k-1}} \!\frac 1{2^{2k-1}}
\sum_{i=1}^{2k}\chi_{1}\dots\hat{\chi_{i}}\dots\chi_{2k}(g)
=\frac 1{2^{2k-1}} \!\sum_{i=0}^{2k-1}(2k-2i) \binom{2k-1}i \!= \!1.
$$
Here $\chi_{1}\dots\hat{\chi_{i}}\dots\chi_{2k}$ denotes the product of
characters
$\chi_{1}\dots{\chi_{i}}\dots\chi_{2k}$ without the character $\chi_{i}.$
Hence
dim$_{\Q}[\Lambda^{2k-1}(\Q^{2k})]^{\Z_2^{\,2k-1}} = 1.$
Finally, one can calculate that  $$\sum_{g\in \Z_2^{\,2k-1}}
\chi_{1}\chi_{2}\dots\chi_{2k}(g) =
\sum_{i=0}^{2k-1}(-1)^{i} \binom{2k-1}i=0.
$$
In the formulas above we use the definition of characters on the
elements of $\Z_2^{\,2k-1}$ considered as elements of $\Z_{2}$-vector space.
\end{proof}

A small modification to the groups $K_n$ defined in
Section~\ref{no} produces examples as required in Remark
\ref{iii}(iii): we replace the $(1,1)$ entries in matrices $C_i$
by $2\times 2$ matrices of the form
$$\begin{array}{rrl}
 & J:=\left[\begin{array}{cc} 0&1\\1&0\end{array} \right]
 & \text{if the entry is a $-1$; and} \\[.4cm]
 \text{the $2\times 2$ identity }
 & \left[\begin{array}{cc} 1&0\\0&1\end{array} \right]
 & \text{if the entry is a 1}.
\end{array}$$
Now the dimension increases by one, while the holonomy group
remains the same. For example, for $n=2$, we had for $K_2$ (the
Klein bottle group), matrix $C_1=\left[\begin{array}{cc}
-1&0\\0&1\end{array} \right]$ and vector $s(C_1)=\frac{e_2}2$; and
now we have a $3\times3$ matrix $\left[\begin{array}{cc}
J&0\\0&1\end{array} \right]$ and vector $\frac{e_3}2$ to obtain a
3-dimensional Bieberbach group which is not of diagonal type.

Hence these are, as claimed, non-diagonal examples in the case
that $k+2$ equals the dimension (where $\Z_2^k$ is the holonomy
group). The verification that these are indeed Bieberbach groups
is easy, and analogous to that for $K_n$. The smaller the $k$
(with fixed dimension), the easier the construction of these
examples.

\section{Outer Automorphism Groups}\label{outaut}

In this section we give an estimate for the order of the outer
automorphism group of GHW Bieberbach groups.
In the case of the groups $K_n$ considered in Section \ref{no}, we prove
that the group
$\,\mathrm{Aut}(K_n)/\mathrm{Inn}(K_n)=: \mathrm{Out}(K_n) \simeq \Z_2^{\,n}.$
We denote by $|A|$ the order of a group $A$.

\begin{theorem}\label{out}
Let $\Gamma$ be a Generalized Hantzsche-Wendt group of dimension $n$.
Then we have $\,|\mathrm{Out}(\Gamma)| \leq 2^{n+1}n!$.
Moreover in the case with non-trivial center we have
$\,|\mathrm{Out}(\Gamma)| \leq 2^{n}(n-1)!.$
\end{theorem}

\begin{proof}
Let $N$ denote the normalizer in $\Gl(n,\Z)$ of the holonomy group
$\Z_2^{\,n-1}$. Then $N$ acts in a natural way on
$H^{2}(\Z_{2}^{\,n-1}, \Z^{n})$.

Let $\alpha \in H^{2}(\Z_2^{\,n-1}, \Z^{n})$ denote
the cohomology class giving rise to the group $\Gamma,$ and
let $N_{\alpha}$ denote its stabilizer in $N$. Then $\Z_2^{\,n-1}$
is a normal subgroup of $N_{\alpha}$
and we have a short exact sequence (cf.\ \cite{Ch}, \cite{HS})
\begin{equation}\label{ast}
0\longrightarrow H^{1}(\Z_2^{\,n-1}, \Z^{n})\longrightarrow
\mathrm{Out}(\Gamma)
\longrightarrow N_{\alpha}/\Z_2^{\,n-1}\longrightarrow 0.
\end{equation}
From \cite{HS} for HW groups, $N = S_{n}\wr \Z_{2}$ is semidirect
product of the symmetric group $S_{n}$ and $\Z_2^{\,n}$. For a
non-oriented GHW group with non-trivial center, $N = \Z_{2}\times
(S_{n-1}\wr \Z_{2})$. In fact, the normalizer $N\simeq
\Z_{2}\oplus N_{1}$ because the image of the holonomy
representation contains a unique element having as its diagonal
entries only one 1 and the remaining are all $-1$'s, (see Lemma
\ref{contains}). Moreover $N_{1}\simeq S_{n-1}\wr \Z_{2}.$ In the
case of non-oriented $\Gamma$ with trivial center we can write the
diagonal group $\Z_2^{\,n}$ in the group $\Gl(n,\Z)$ as a sum of a
holonomy group $\Z_2^{\,n-1}$ and a set $(-Id)\Z_2^{\,n-1}$. Hence
$N_{\Gl(n,\Z)}(\Z_2^{\,n-1}) \subset N_{\Gl(n,\Z)}(\Z_2^{\,n})
\simeq S_{n}\wr \Z_{2}$. From (\ref{ast}) we only have to prove
that $| H^{1}(\Z_2^{\,n-1},\Z^{n}) |$ equals $2^{n}$ in the case
of trivial center and $2^{n-1}$ in case of non-trivial center. But
this is proved in \cite{Ch}.
\end{proof}

\begin{cor}
For $n\geq 2$, $\mathrm{Out}(K_n) \simeq \Z_2^{\,n}.$
\end{cor}

\begin{proof} From above the normalizer $N$ of the holonomy representation
is equal to $\Z_{2}\times S_{n-1}\wr \Z_{2},$ where the
permutation group acts on the first $(n-1)$ coordinates. Let $s\in
H^{2}(\Z_2^{\,n-1}, \Z^{n})$ be a cocycle which corresponds to the
group $K_n$ of Section~\ref{no}. Then $N_{s} = \Z_2^{\,n}.$ In
fact it is enough to observe that for any non-trivial $\sigma \in
S_{n-1}$, \,$\sigma s \neq s.$ We have
$$(\sigma s)(B_{i}) = \sigma s(\sigma^{-1} B_{i}\sigma) = \sigma
s(B_{\sigma^{-1}(i)}),$$ for any generator $B_{i}\in \Z_2^{\,n-1}$
and $i= 1,2,..., n-1,$ see \cite{HS}. We are looking for the
permutation $\sigma \in S_{n-1}$ such that the cocycle $\sigma s -
s$ is trivial. We have two cases. First assume that
$\sigma^{-1}(i) + 1 < n.$ Hence $\sigma^{-1}(i) + 1 =
\sigma^{-1}(i + 1).$ From  $\sigma^{-1}(i) + 1 = n$ it follows
that $i = n-1.$ Summing up we have $\sigma \in N_{s}$ if and only
if, for $i< (n-1)$, \,$\sigma(i + 1) = \sigma(i) + 1$ and
$\sigma(n-1) = n-1.$ Hence for $i = n-2$ we have $\sigma(n-2) =
n-2$ and by induction $\sigma$ is identity. To finish the proof we
observe that $N_{s}/\Z_2^{\,n-1} = \Z_{2}$ and it is generated by
the image of $-Id$. Now it is enough to apply the short exact
sequence (\ref{ast}).
\end{proof}

\begin{remark} {\rm From \cite{Ch}, for GHW manifolds $M^{n}$ of dimension
$n$ with first Betti number $\beta_{1}$, the group
Aff$(M^{n})/(S^{1})^{\beta_{1}}$ of affine diffeomorphisms of $M^{n}$ is
isomorphic to the group $\mathrm{Out}(\pi_{1}(M^{n})).$
}\end{remark}

\section{A Graph Connecting the GHW Manifolds}\label{graph}

In this section we start the discussion of relations between the generalized
Hantzsche-Wendt manifolds of different dimensions.
For manifolds with first Betti number one we have relations
from E.\ Calabi's construction (see Proposition \ref{chain} and \cite{Wo}, p.\
125).
But we have also a lot of links between GHW groups of different dimensions
in the case when one of them has trivial center. We must say that
this subject needs further research and we want to come back to it
in the future.
Most of the observations in this section are corollaries from previous ones,
and we shall omit some proofs.

Let us start defining a graph {\bf G} as follows:
\begin{itemize}
\item The vertices of {\bf G} are the GHW groups (or manifolds).
\item We say that a vertex has dimension $n$ if it corresponds to a GHW
manifold
of dimension $n$.
\item Edges are defined only between vertices in consecutive dimensions.
 \item There is an edge between two vertices (in consecutive dimensions) if one
of the corresponding GHW groups
is a subgroup of the other (in an abstract sense, i.e. there is a
monomorphism from one into the other).
\end{itemize}

From Corollary \ref{subgroup} we see that in {\bf G} there is
always (at least) an edge `down' from each vertex. As a direct
consequence we have

\begin{cor}
(i) {\it The fundamental group of the Klein bottle is a subgroup of every GHW
group.} \\
(ii) {\it The graph ${\bf G}$ is connected.}
\end{cor}

\begin{figure}[htb]
\input{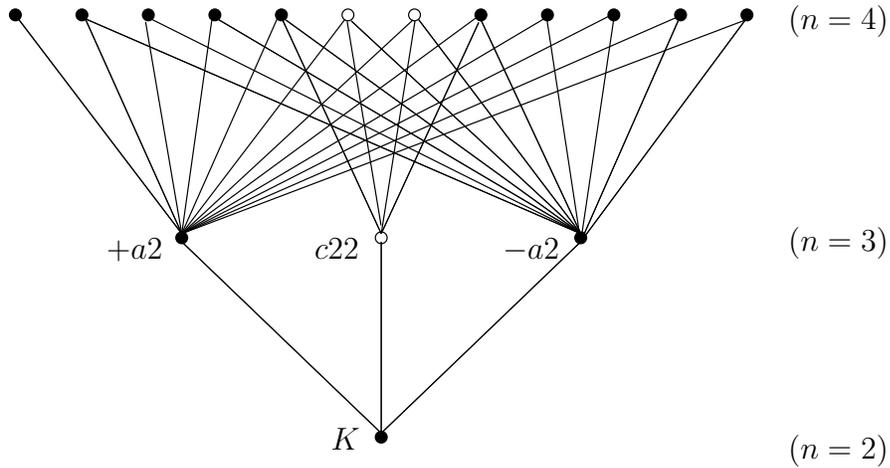}
\caption{The graph of GHW manifolds up to $n=4$.}
\label{fig}
\end{figure}

Moreover from Figure \ref{fig} we see that {\bf G} is not a tree.

On the other hand, from the method in the proof of Theorem \ref{exist} it is
possible to
deduce that any GHW group in dimension $n$ is subgroup of a GHW group in
dimension $n+1$.
Moreover, we can prove that it is indeed subgroup of at least two different
GHW groups in dimension $n+1$ (we omit the proof for brevity).
In the graph ${\bf G}$ this means that from any vertex there are at least two
edges `up'.

In Figure \ref{fig}, a point $\circ$ (resp.\ a bullet $\bullet$) denotes a GHW
manifold with first
Betti number zero (resp.\ one);
$K$ is the Klein bottle;
$+a2$ and $-a2$ are the {\it first} and {\it second amphidicosm}
respectively,
they are the two non-oriented GHW manifolds of dimension 3;
$c22$ is the {\it didicosm} (or the oriented Hantzsche-Wendt 3-dimensional flat
manifold).
The new names are due to a proposal by John H.\ Conway to name the ten flat
3-manifolds according to their properties, which we adopt (this
terminology is introduced in \cite{CR}).

Other edges in Figure \ref{fig} are obtained by extensions to the methods of
theorems \ref{diag} and \ref{exist}.

\vskip.03cm

There are many other properties of ${\bf G}$ to be explored.
For instance, a {\it distance} between vertices  can be defined in the standard
way for graphs.

In dimension 5 there are only two HW manifolds (cf.\ \cite{MRc}).
With the method of Theorem \ref{diag}, it is not difficult to
prove that they have both $\G_4$ (defined in Section
\ref{examples}) as a subgroup, and therefore the distance between
them is just two. We also have

\begin{prop}
The fundamental group of the didicosm (or Hantzsche-Wendt 3-manifold) is
a subgroup of every GHW group with trivial center.
\end{prop}
\begin{proof}
Denote by $\G$ the given GHW group. By Lemma \ref{contains} there
is $(B,b)\in\G$ such that (after permuting coordinates)
$B=\text{diag}(1,-1,\dots,-1)$ and $b=(\frac 12,*,\dots,*)$, where
each $*$ is in $\{0,\frac 12\}$. Since $\G$ has trivial center
there is another element $(C,c)\in\G$ with
$C=\text{diag}(-1,1,-1,\pm1,\dots,\pm1)$ and the second coordinate
of $c$ equals $1/2$. Now the integral representation corresponding
to $\langle B, C, BC \rangle$ is the direct sum of the three
non-trivial characters of $\Z_2\oplus\Z_2$ (with certain
multiplicities). Then it is not difficult to check that the
subgroup of $\G$ generated by $(B,b)$ and $(C,c)$ is isomorphic to
the fundamental group of the didicosm.
\end{proof}

In \cite{GS}, the analogous result was proved for Bieberbach groups with
holonomy
group $\Z_2\oplus\Z_2$ (instead of GHW groups) and trivial center (see also
\cite{Ti}).

\vskip0.3cm

By (\ref{lee}) the groups $K_n$, $n\ge 2$, form an infinite chain in {\bf G}.
Also the groups $\G_n$, $n\in \N$, $n\ge 3$ from Section \ref{examples} form an
infinite chain by Proposition \ref{mono}.

Finally, we remark that the edges in {\bf G} sometimes correspond to normal
subgroups and
sometimes not (or they can correspond to both simultaneously). For instance,
the
chain obtained
by Proposition \ref{mono} is not normal.
Also, $\pi_1(K)$
can be injected as a normal and as a non-normal subgroup in $\pi_1(+a2)$ and in
$\pi_1(-a2)$ but
it is never a normal subgroup of $\pi_1(c22)$.


\noindent{\bf Some problems and questions:}

\begin{itemize}
\item[{\bf 1.}] Are there relations between Fibonacci groups and
the fundamental groups of Hantzsche-Wendt manifolds? (cf.\ \cite{Sz3}) 
\item[{\bf 2.}] Are there other geometric structures
on oriented GHW manifolds?
Is it possible to introduce some invariants on them? (in analogy with
three dimensional rational homology sphere as Casson-Walker and
Seiberg-Witten invariants) 
\item[{\bf 3.}] Is it possible to
improve the upper bound given in Theorem~\ref{out}? Could the
factor $n!$ be replaced by $n$? 
\item[{\bf 4.}] What are the
geometrical and topological properties of the `graph' of GHW
manifolds?
\end{itemize}

\end{document}